\documentclass{article}[12pt]
\usepackage{graphicx}
\usepackage{amssymb}
\usepackage{amsthm}
\usepackage{mathrsfs}
\usepackage{wrapfig}
\usepackage{tikz} 
\usepackage{tikz-3dplot}
\usepackage{bm}
\usepackage{float} 
\usepackage{amsmath}
\newtheorem{theorem}{Theorem}
\newtheorem{corollary}{Corollary}
\newtheorem{lemma}{Lemma}

\newtheorem{definition}{Definition}

\begin{document}
\title{
Common Limiting Structure Between Classical Wave Theory 
 and Renormalization Theory in Multiple-Scattering Problems}

\author{Terumi Touhei 
\thanks{Tokyo University of Science,\ \ email:\    ttouhei@rs.tus.ac.jp}}
\date { }
\maketitle

\begin{abst}
{}\\
In this study, we investigate Dirac-type scatterers from the viewpoints of
classical wave theory and renormalization theory.
We show that the diagonal self-interaction contribution can be represented
asymptotically by the external field evaluated at the scatterer center as the
spatial scale of the scatterers tends to zero, provided that the transmission
and emission coefficients are properly defined.
Through this analysis, the bare coupling constant that appears in the governing
equation is naturally redefined in terms of the scattering strength, from which
the Foldy equation emerges as the leading-order interaction equation within the
framework of classical wave theory. This redefinition is not apparent from the
governing equation alone.
We further identify two distinct types of strong resolvent convergence for the
resolvent constructed from the Foldy equation. One corresponds to convergence
to the free Hamiltonian, indicating that the interaction effects vanish as the
spatial scale of the scatterers tends to zero. The other corresponds to
convergence to a nontrivial Hamiltonian for which the interaction effects
remain in the zero-size limit.
Although the limiting Hamiltonian in renormalization theory
is formally described by a family of self-adjoint extensions
with an arbitrary real parameter,
the limiting procedure selects a specific value from this family.
This phenomenon can also be observed in 
classical wave theory using the Foldy resolvent.
Therefore,  the limiting Foldy resolvents admit a natural interpretation
within the framework of renormalization theory. These results reveal a common
limiting structure shared by classical wave theory and renormalization theory
in multiple-scattering problems.

\vspace{0.2cm}

\noindent
{\bf Key words}: renormalization theory,
classical wave theory,
Foldy algebraic equation,
Dirac-type scatterers,
scattering strength,
strong resolvent convergence,
multiple scattering
\end{abst}
\section{Introduction}
The analysis of multiple-scattering problems is essential in a wide range of fields,
including non-destructive testing, geophysical exploration, wave propagation in
heterogeneous media, and medical imaging. Such problems arise not only in
classical wave theory but also in quantum mechanics.

For scatterers whose spatial scale is much smaller than the wavelength,
Dirac-type scatterers (point interaction models) provide an important idealized
description in both classical and quantum settings. On the other hand,
less singular models with finite spatial extent are often employed as more
realistic representations of physical scatterers.

For such finite-size scatterers, Foldy (1945) proposed an algebraic equation that describes the interactions among scatterers through
the concept of an external field and developed a statistical framework for
randomly distributed scatterers. These foundational ideas have inspired a large
body of subsequent research in classical wave theory.

For example, Lax (1951, 1952) extended Foldy's framework to elastic and
inelastic scattering. Waterman and Truell (1961) investigated the validity of
the corresponding integral equation using a statistical approach. Angel and
Koba (1998) studied scattering in a cracked slab, while Aristégui and Angel
(2002) examined the problem of isotropic point scatterers. More recently,
Challa and Sini (2014) provided a mathematical justification of the Foldy-Lax
approximation for many small rigid obstacles of arbitrary shape.

Although Foldy's framework is widely used for finite-size scatterers, its relation
 to Dirac-type scatterers and renormalization theory remains less clear.
In Martin's monograph (2006), the Foldy equation is not applied directly
to Dirac-type scatterers. Instead, alternative methods, such as the localized
perturbation method (Lax, 1954), are surveyed to address divergences caused by self-interaction
in point-scatterer models (see Martin (2006), Chapter 8).

In contrast, several studies have applied the Foldy equation to
Dirac-type scatterers. For example, Hu and Sini (2013) treated the Foldy equation
as an approximation of the wavefield in the presence of Dirac-type scatterers and
justified its use for two- and three-dimensional elastic wavefields using a renormalization
framework from quantum mechanics (Albeverio et al., 1988). 
In a follow-up study, Hu et al. (2020) extended
renormalization techniques to point-like and extended obstacles in 
elastic scattering problems.

The present study is motivated by the renormalization framework developed by
Albeverio et al. (1988).  According to their monograph, the resolvent of a Hamiltonian
 formally containing Dirac-type interactions cannot be directly defined, 
 and a renormalization procedure is required to construct the 
 corresponding self-adjoint Hamiltonian.
This naturally raises the question of how such observations affect the applicability of the
Foldy algebraic equation in classical wave theory.

In this study, we begin with the equation governing sound-wave propagation
in the presence of \(N\) small spherical scatterers, 
where the effects of the gradient of the mass density of the medium
are negligible (Colton and Kress, p.~267).  
Despite this simplification, this equation has been used in the analysis of 
non-invasive organ examination (see, for example, Ueda, 2024).  

We investigate the scattered wavefield in the limit \(a \to 0\) using classical
wave theory, where \(a\) denotes the spatial scale of the scatterers.
The governing equation can also be interpreted as a Schr\"odinger equation with
\(N\) compactly supported potential wells.
This observation makes it possible to compare classical wave theory 
with the renormalization technique developed in quantum mechanics.

We show that the diagonal (self-interaction) contribution can be represented 
asymptotically by the external field evaluated at the scatterer center, provided
that the transmission and emission coefficients are properly defined. 
During the limiting process \(a \to 0\), the bare coupling constant that 
appears in the original 
governing equation must be redefined in terms of the resulting 
scattering strength.
Consequently, the Foldy equation emerges naturally as the leading-order 
interaction equation within the classical wave theory.  

The most important finding of this study is that the classical wave theory
 and the renormalization framework exhibit a common limiting structure.
This observation helps bridge the conceptual gap between classical wave theory
and quantum mechanical treatments based on renormalization.

\section{Basic Equations and Problem Definition}

In this study, we discuss multiple-scattering problems involving
very small spherical objects, where the size of each scatterer is much 
smaller than the wavelength of the background wavefield and the 
scatterers are mutually disjoint. 
We consider the following type of equation:
\begin{eqnarray} 
\Bigl(
-\Delta - \sum_{\alpha=1}^{N}
 \mu_{\alpha}  \frac{X_a (\bm{r}-\bm{r}_{\alpha} )}{a^d}
-k^2 \Bigr) u(\bm{r}) &= & \delta ( \bm{r}-\bm{r}_F),
\nonumber \\
 (\bm{r}\in \mathbb{R}^3,\ a>0,\ d=1,2,3),& &
\label{eqn:e1}
\end{eqnarray}
where
\(\Delta\) is the Laplacian, \(\delta(\cdot)\) is the Dirac delta function,
\(\bm{r}_F\) is the location of the point source generating the incident wavefield,
and \(X_a (\cdot)\) is the indicator function defined by
\begin{equation}
  X_a (\bm{r}- \bm{r}_{\alpha})=
  \begin{cases}
    1 & \text{if } |\bm{r} -\bm{r}_{\alpha}| \le a, \\
    0 & \text{otherwise}.
  \end{cases}
\end{equation}  
Here, \(a\) is the radius of the spherical 
scatterer, \(\bm{r}_{\alpha}\) denotes the center of 
the scatterer, and \(d\) is a parameter characterizing the singularity of the scatterer.

In the classical mechanics regime, Eq.~(\ref{eqn:e1}) describes the scattering of sound waves
due to the presence of spherical scatterers, where
the effects of the gradient of the mass density of the medium
are negligible (Colton and Kress, p.~267).  

Equation~(\ref{eqn:e1}) can also be understood as a Schr{\"o}dinger 
equation with \(N\) potential wells.
In quantum mechanics, the equation
\begin{equation}
\Bigl(
-\Delta - \sum_{\alpha=1}^{N}
\mu_{\alpha} \delta (\bm{r}-\bm{r}_{\alpha})
-k^2 \Bigr) u(\bm{r}) = \delta (\bm{r} - \bm{r}_F )  
\label{eqn:e2}
\end{equation}
has been studied in detail (see, for example, Albeverio et al., 1988), where 
\(\mu_{\alpha}\) is referred to as the coupling constant. 
In the context of  quantum mechanics,
Eq.~(\ref{eqn:e1}) describes a more realistic 
and less singular interaction model than Eq.~(\ref{eqn:e2}).

In this study, we analyze Eq.~(\ref{eqn:e1})
using classical wave theory,
including the investigation of 
the asymptotic properties of the wavefield and expansions in
spherical harmonics. 
We compare the results obtained from the classical
wave theory  with those obtained from quantum mechanics
discussed in the literature (Albeverio et al., 1988).  
Note that Albeverio et al. (1988) considered only the case \(d=2\), which 
is associated with the unitary scaling group. In this study, 
we extend the analysis to the cases \(d=1\) and \(d=3\).

For simplicity, we assume that the coupling constants are real and positive.  
Green's function for the background medium is given by
\begin{equation}
G(\bm{r}, \bm{r}') = \frac{1}{4 \pi |\bm{r} - \bm{r}'|}
\exp\left( -i k |\bm{r} - \bm{r}'| \right), \qquad 
(\bm{r}, \bm{r}' \in \mathbb{R}^3).
\end{equation}

In the remainder of this study, we first apply classical 
wave theory approaches to Eq.~(\ref{eqn:e1}) 
and show that the diagonal (self-interaction) term is  
asymptotically equivalent to the external field evaluated at the scatterer center.
As a result,  the Foldy equation is naturally derived as the 
interaction equation within the classical framework, provided that 
the emission and transmission coefficients for each scatterer are properly defined. 
Finally, we compare these results with
those obtained from renormalization theory.
We show that, as the spatial scale of the scatterers tends to zero, the
Foldy algebraic equation and the renormalization technique can be understood 
within a common limiting structure.

\section{Classical Wave Theory for Finite Spatial Scale Scatterers}

\subsection{Possibility of Monopole Wavefield Approximation}

At this stage, within the framework of classical wave theory, 
we focus on the \(\alpha\)-th scatterer 
in a multiple-scattering system consisting of \(N\) scatterers. 
We notice that there are three types of waves with respect to the 
\(\alpha\)-th scatterer, namely those propagating toward it, those 
transmitted into it, and those scattered from it. 
They are denoted by \(u_{\alpha}(\bm{r})\), \(u_{\alpha}^{\downarrow}(\bm{r})\), and
\(u_{\alpha}^{\uparrow}(\bm{r})\), respectively. 

The waves propagating to  the \(\alpha\)-th scatterer consist of the sum of the 
scattered waves emitted from all other scatterers and the incident wave. 
This combined field is referred to as the \textit{external field} (Foldy, 1945) for 
the \(\alpha\)-th scatterer.   
We introduce a local coordinate system centered at the scatterer, 
described in terms of the spherical coordinates \((r, \theta, \varphi)\).

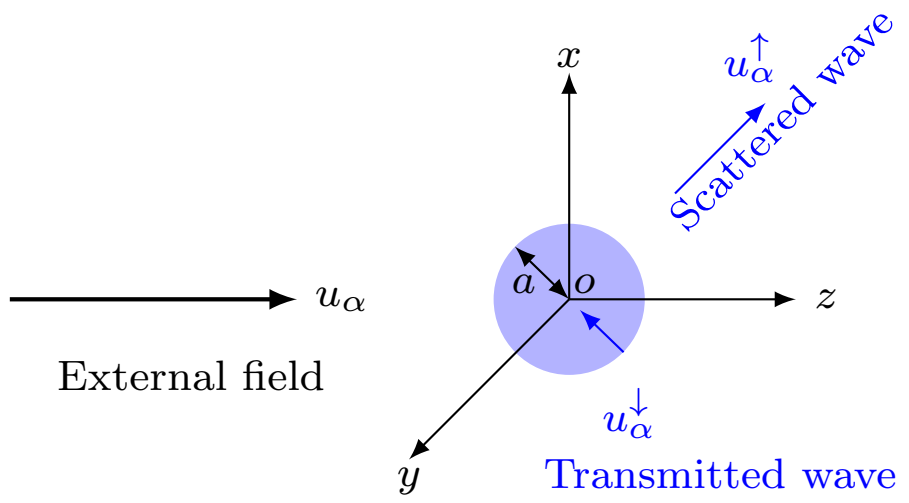
\begin{figure}[h]
\begin{center}
\begin{tikzpicture}[inner sep=0pt]

\node at (0.5,-0.5){\footnotesize {External field}};
\node at (1.5,0.0){\footnotesize {\( u_{\alpha} \)}};

\draw[->, thick, >=latex] (-0.7,0.0)--(1.2,0.0);
  \fill [blue, opacity=0.3](3,0) circle[radius=0.5];
\draw [->,>=latex]( 3,0)--(1.94,-1.06);
\draw [->,>=latex]( 3,0)--(3,1.5);
\draw [->,>=latex]( 3,0)--(4.5,0);

\node at (3.10,0.1){\footnotesize $o$};
\node at(3,1.6){\footnotesize $x$};
\node at (1.94,-1.2){\footnotesize $y$};
\node at (4.7,0){\footnotesize $z$};
\draw [<->,black, >=latex] (3,0) --(2.64,0.35);
\draw [->,blue, >=latex] (3.36,-0.35 ) --(3.07,-0.07);

\node at (3.4,-0.75) {\color{blue} 
\footnotesize{\( u_{\alpha}^{\downarrow}\)}};
\node at (4.0,-1.15){\color{blue}
\footnotesize{Transmitted wave}};
\node at (2.7,0.1){\footnotesize \color{black}$a$};

\draw [->, blue, >=latex](3.7,0.7)--(4.3,1.3);
\node at (4.4,1.2)[rotate=45] {\color{blue} \footnotesize{Scattered wave}};
\node at (4.2,1.6) {\color{blue} 
\footnotesize{\( u_{\alpha}^{\uparrow}\)}};
\end{tikzpicture}

\caption{
  Local coordinate system for \(\alpha\)-th scatterer and
  definitions of \(u_{\alpha}\), \(u_{\alpha}^{\uparrow} \), and
  \(u_{\alpha}^{\downarrow}\).
 The scatterer's center is at the origin of the local coordinate system.  
\label{fig:concept}}
\end{center}
\end{figure} 

We express the spherical harmonics expansions of 
\(u_{\alpha} (\bm{r} )\), \(u_{\alpha}^{\downarrow}( \bm{r} )\), and
\(u_{\alpha}^{\uparrow} (\bm{r} )\) as
\begin{align}
u_{\alpha} ( \bm{r} ) &=
  \sum_{n=0}^{\infty} \sum_{|m| \le n} F_{\alpha n}^{\ m}
Y_{n}^{m}(\theta, \varphi)\,
j_{n} (kr),
\label{eqn:Eextf}
\\
u_{\alpha}^{\downarrow} ( \bm{r} ) &=
  \sum_{n=0}^{\infty} \sum_{|m| \le n} A_{\alpha n}^{\ m}
Y_{n}^{m}(\theta, \varphi)\,
j_{n} (k_{\alpha} r),
\label{eqn:Intf}
\\
u_{\alpha}^{\uparrow} ( \bm{r} ) &=
  \sum_{n=0}^{\infty} \sum_{|m| \le n} B_{\alpha n}^{\ m}
Y_{n}^{m}(\theta, \varphi)\,
h_{n}^{(2)} (k r),
\label{eqn:Emif}
\end{align}
where \(j_n\) is the spherical Bessel function of the first kind, 
\(h_n^{(2)}\) is the spherical Hankel function of the second kind, and
\(Y_n^m\) is defined by
\begin{equation}
Y_{n}^{m}(\theta, \varphi)=
\sqrt{\frac{ (2n+1)(n-m)!}{4 \pi (n+m)!}}
P_n^m(\cos \theta)\, e^{i m \varphi},
\label{eqn:surf}
\end{equation}
where \(P_n^m\) represents the associated Legendre polynomials. 
The wavenumber \(k_{\alpha}\) is defined by
\begin{equation}
k_{\alpha}^2 = k^2 + \frac{\mu_{\alpha}}{a^d}.
\label{eqn:def_kalpha}
\end{equation}

The coefficients \(A_{\alpha n}^{\ m}\) and \(B_{\alpha n}^{\ m}\) are
 determined from \(F_{\alpha n}^{\ m}\)
by imposing the continuity of the wavefield and its radial derivative at the interface \(|\bm{r}|=a\):
\begin{align}
\lim_{|\bm{r}| \to a+0} \bigl(u_{\alpha} (\bm{r}) +
u_{\alpha}^{\uparrow} (\bm{r})\bigr)
&= \lim_{|\bm{r}| \to a-0} u_{\alpha}^{\downarrow} ( \bm{r} ),
\label{eqn:ub_aE1}
\\
\lim_{|\bm{r}| \to a+0} 
\frac{\partial}{\partial r} \bigl(u_{\alpha} (\bm{r}) + u_{\alpha}^{\uparrow} (\bm{r})\bigr)
&= \lim_{|\bm{r}| \to a-0} 
\frac{\partial}{\partial r} u_{\alpha}^{\downarrow} ( \bm{r} ).
\label{eqn:ub_aE}
\end{align}

We now examine the wavefield under the condition \(|ka|\ll 1\) and obtain the following lemma.

\begin{lemma}
The scattered wavefield is dominated by the \(n=0\) term when \(|ka| \ll 1\).
Therefore, the monopole wavefield approximation is possible.
\end{lemma}

\noindent
\begin{proof}
According to Eqs. (\ref{eqn:ub_aE1}) and   (\ref{eqn:ub_aE}), 
we can determine the coefficients $A_{\alpha n}^{\ m}$ and 
$B_{\alpha n}^{\ m}$ as
\begin{equation}
\left (
\begin{array}{c}
A_{\alpha n}^{\ m} \\
B_{\alpha n}^{\ m}
\end{array} \right )
= \frac{1}{\Delta_n (a)}
\left (
\begin{array}{cc}
\partial_r h_n^{(2)} (k a ) & -h_n^{(2)} (k a ) 
\\
\partial_r j_n (k_\alpha a )  & - j_n (k_\alpha a ) 
\end{array} \right )
\:
\left ( \begin{array}{c}
j_n (k a ) \\
\partial_r j_n (k a) 
\end{array} \right ) F_{\alpha n}^{\ m},
\label{eqn:sol_AB}
\end{equation}
where 
\begin{equation}
\Delta_n (a)
= j_n (k_\alpha a) \partial_r h_n^{(2)} (k a) -
h_n^{(2)} (k a) \partial_r j_n (k_\alpha a).
\label{eqn:deltan}
\end{equation}
Here, we employ the following notation for the derivative of
the spherical Bessel functions: 
\begin{equation}
\partial_r j_n (ka )= \Bigl [\frac{\partial j_{n} (kr)}{\partial r} 
\Bigr ]_{r=a}
\end{equation}

We introduce the following asymptotic properties 
of the spherical Bessel functions:
\begin{eqnarray}
  j_{n}(ka)& \sim & \frac{ (ka)^{n}}{(2n+1)!!} ,\ \ (ka \to 0) 
\nonumber \\
  h_n^{(2)} (ka) & \sim & \frac{i (2n-1) !!}{(ka)^{n+1}} ,\ \ ( ka \to 0)
\nonumber \\
  j_{n} ( z) & \sim & \frac{1}{z}\cos 
  \Bigl ( z-\frac{(n+1)\pi}{2} \Bigr ),
  \ \ ( |z| \to \infty)
\label{eqn:e18}
\end{eqnarray}
as well as the definition of $k_{\alpha}$ represented by Eq. (\ref{eqn:def_kalpha})
into Eq. (\ref{eqn:sol_AB}). 
Since the emission of the scattered waves is determined by the coefficient
\( B_{\alpha n}^{\ m} \) , 
we only have to present the following asymptotic evaluation of 
\( B_{\alpha n}^{\ m} \) :
\begin{equation}
B_{\alpha n}^{\ m} =
\begin{cases}
O((ka)^2), & (d=1, n=0), \\
O((ka)^{2n+1}), & (d=1, n \ge 1), \\
O((ka)^{2n+1}), & (d=2,3, n \ge 0).
\end{cases}
\end{equation}
Therefore, the scattered wavefield can be approximated by the \(n=0\) 
term alone, which completes the proof.
\end{proof}

\vspace * {0.1cm}                          

\noindent
\textbf{Remark 1.}
Instead of Eq.~(\ref{eqn:ub_aE}), if we employ
\begin{eqnarray*} 
\lim_{|\bm{r}| \to a+0}  
\frac{1}{\rho_e} \frac{\partial}{\partial r} \bigl(u_{\alpha}
(\bm{r}) + u_{\alpha}^{\uparrow} (\bm{r})\bigr)
& = & \lim_{|\bm{r}| \to a-0}
\frac{1}{\rho_i}
\frac{\partial}{\partial r} u_{\alpha }^{\downarrow} ( \bm{r} ),
\nonumber \\
\qquad (\rho_e \neq \rho_i),
\end{eqnarray*} 
where \(\rho_e\) and \(\rho_i\) are the mass densities
outside and inside the scatterer, respectively,
then a dipole field (i.e., the \(n=1\) term) is induced
in the transmission problem (Martin, 2006, p.~302). 

However, if we restrict ourselves to the case \(\rho_e = \rho_i\),
the above lemma ensures that 
only the monopole field (\(n=0\)) needs to be considered when 
\(|k a| \ll 1\).
The condition \(\rho_e = \rho_i\) is consistent with the assumption
that the gradient of the mass density is negligible, as stated previously.
Furthermore, this assumption leads to a formulation analogous to
the Schr{\"o}dinger equation with potential wells
and enables comparison with results
from quantum mechanics. 

\subsection{Transmission and Emission Coefficients }   

Under the monopole wavefield approximation,
let \(A_{\alpha}\), \(B_{\alpha}\), and \(F_{\alpha}\) denote
\(A_{\alpha 0}^{\ 0}\), \(B_{\alpha 0}^{\ 0}\), and \(F_{\alpha 0}^{\ 0}\),
respectively, in Eqs.~(\ref{eqn:Eextf})--(\ref{eqn:Emif}), which are
reduced to
\begin{eqnarray}
u_{\alpha}^{\uparrow} ( \bm{r} ) &=&
B_{\alpha}\,
Y_{0}^{0}(\theta, \varphi)\,
h_{0}^{(2)} (k r),
\label{eqn:Emif1}
\\                                      
u_{\alpha}^{\downarrow} ( \bm{r} ) &=&
A_{\alpha}\,
Y_{0}^{0}(\theta, \varphi)\,
j_{0} (k_{\alpha} r),
\label{eqn:Intf1}   
\\
u_{\alpha} (\bm{r} ) &=&
F_{\alpha}\,
Y_{0}^{0}(\theta, \varphi)\,   
j_0 (kr).
\label{eqn:Exterf1}
\end{eqnarray}

We refer to \(A_{\alpha}\) and \(B_{\alpha}\) as
the {\it transmission and  emission coefficients}
with respect to the \(\alpha\)-th scatterer, respectively. 
Our task at this point is to evaluate the asymptotic properties of
the transmission and emission coefficients when \(|ka| \ll 1\). 

In the following, we consider four cases according to the parameter \(d\) and
the coupling constant \(\mu_{\alpha}\), as summarized in 
Table~\ref{tab:tb1}.
In particular, for the case \(d=2\), we distinguish two subcases depending on whether 
\(\cos(\sqrt{\mu_{\alpha}})=0\) or \(\cos(\sqrt{\mu_{\alpha}})\neq 0\), 
as will become clear from the asymptotic behavior of 
\(A_{\alpha}\) and \(B_{\alpha}\).

\begin{table}[H]
  \caption{Cases Considered in Analysis.
    \label{tab:tb1}
  }
  \begin{center}
  \begin{tabular}{c|l}
    \hline
    Case & Conditions\\
    \hline
    1    & \(d=1\)\\
    2a   &  \(d=2\) ,\ \( \cos (\sqrt{\mu_{\alpha}} )\ne 0 
    \)\\
    2b 	&  \(d=2\) ,\ \( \cos (\sqrt{\mu_{\alpha}} ) =0    \)\\ 
    3   &  \(d=3\)\\
    \hline
  \end{tabular}
  \end{center}
  \end{table}

The expressions for \(A_{\alpha}\) and \(B_{\alpha}\) for cases~1, 2a, and~3 are given by
\begin{equation}
  \left \{
  \begin{array}{rcl}
    A_{\alpha} &= & {\displaystyle
      \frac{1}{\mathcal{N}_{\alpha}} F_{\alpha}}, \\
    B_{\alpha} &= & {\displaystyle (-ika)\frac{\mathcal{M}_{\alpha}}
    {\mathcal{N}_{\alpha}} F_{\alpha}},\ \ \ 
  \end{array} \right .
\label{eqn:defABalpha}
\end{equation} 
where $\mathcal{N}_{\alpha}$ and $\mathcal{M}_{\alpha}$ are given by
\begin{eqnarray}
  \mathcal{N}_{\alpha}
  & = & \left \{
\begin{array}{ll}
  1 & {\rm{in\ \ case\ 1}}\\
\cos (\sqrt{\mu_{\alpha} }) 
& {\rm{in\ \ case\ 2a}} \\
  \cos ( \sqrt{\mu_{\alpha}/a }) 
& {} \\
 +ika \Bigl (j_0 \bigl ( \sqrt{\mu_{\alpha} /a} \bigr ) - 
\cos ( \sqrt{\mu_{\alpha}/a }   ) \Bigr )
& {\rm{in\ \ case\ 3}}
\end{array} \right.
\nonumber \\
& & 
\label{eqn:eq21a} \\
  \mathcal{M}_{\alpha}
  & = & \left \{
\begin{array}{ll}
\mu_{\alpha} a/3  & \hspace * {0.95cm} {\rm{in\ \ case\ 1}}\\
\sqrt{\mu_{\alpha}} j_1  (\sqrt{\mu_{\alpha}} )
& \hspace * {0.95cm} {\rm{in\ \ case\ 2a}} \\ 
 j_0 \bigl ( \sqrt{\mu_{\alpha} /a} \bigr ) - 
\cos \left( \sqrt{\mu_{\alpha}/a }  \right ) 
& \hspace *{0.95cm}  {\rm{in\ \ case\ 3}}
\end{array} \right .
\nonumber \\
& & 
\label{eqn:eq22a}
\end{eqnarray}

For case 2b, \(A_{\alpha}\) and \(B_{\alpha}\) are given by
\begin{equation}
  \left \{
  \begin{array}{rcl}
    A_{\alpha} &= & -\displaystyle{
      \frac{i}{k j_0 (\sqrt{\mu_{\alpha}})     }
   \Bigl (
      \frac{1}{a }
      +\frac{\mathcal{C}_0}
      {  j_0 (\sqrt{\mu_{\alpha} })} 
        \Bigr ) F_{\alpha}},\\
    B_{\alpha} &= & -\displaystyle {\Bigl (1+\frac{a\: \mathcal{C}_0 }
      {
          j_0 (\sqrt{\mu_{\alpha}})} \Bigr )} F_{\alpha},
  \end{array} \right .
\label{eqn:defABalpha2}
\end{equation} 
where
\begin{equation}
  \mathcal{C}_0=
  \frac{ik}{2}
  \Bigl (
  \frac{ j_0{'}
    ( \sqrt{\mu_{\alpha}} )
    -\sqrt{\mu_{\alpha}} j_1{'}
    ( \sqrt{\mu_{\alpha}} )}{
    \sqrt{\mu_{\alpha}} } + j_0 
    ( \sqrt{\mu_{\alpha}} ) \Bigr ).
\end{equation}

Note that Eq.~(\ref{eqn:defABalpha}) corresponds to the leading-order term 
in the asymptotic expansion of Eq.~(\ref{eqn:sol_AB}) and
 Eq.~(\ref{eqn:defABalpha2}) corresponds to the next-order term 
 under the condition \(\cos(\sqrt{\mu_{\alpha}})=0\).

\subsection{Definition of Scattering Strength}

Now, let us move away from the local coordinate system
and recall that \(\bm{r}_{\alpha} \) is the center of the \(\alpha\)-th scatterer. 
Then, Eqs.~(\ref{eqn:Emif1}) and (\ref{eqn:Intf1})
can be expressed as
\begin{eqnarray} 
 u_{\alpha}^{\uparrow} (\bm{r} ) &=& \frac{1}{\sqrt{4\pi}} 
 h_0^{(2)} \bigl ( k  | \bm{r} -\bm{r}_{\alpha} | \bigr ) B_{\alpha}, \quad
\bigl (  | \bm{r} -\bm{r}_{\alpha} |>a \bigr ),
\label{eqn:Emif2} \\
u_{\alpha}^{\downarrow} (\bm{r} ) &=& \frac{1}{\sqrt{4\pi}} 
 j_0 \bigl ( k_{\alpha}  | \bm{r} -\bm{r}_{\alpha} | \bigr ) A_{\alpha}, \quad
\bigl (  | \bm{r} -\bm{r}_{\alpha} |<a \bigr ).
\label{eqn:Intf2}
\end{eqnarray}
According to Eq.~(\ref{eqn:Exterf1}), the external field \(u_{\alpha}\) at \(\bm{r}=\bm{r}_{\alpha} \) is expressed as
\begin{eqnarray}
u_{\alpha} ( \bm{r}_{\alpha} )&=&F_{\alpha} Y_{0}^{0} (\theta, \varphi)
\left [  j_{0} (kr )\right ]_{r=0} \nonumber \\
&=& \frac{1}{\sqrt{4\pi}} F_{\alpha}.
\label{eqn:Exterf2}
\end{eqnarray}

Combining the identity
\begin{equation}
-ik\, h_{0}^{(2)} \bigl(k|\bm{r} -\bm{r}_{\alpha} | \bigr)
= \frac{e^{-ik|\bm{r} - \bm{r}_{\alpha}|}}
{|\bm{r}- \bm{r}_{\alpha}|}
\label{eqn:hn25}
\end{equation}
with Eqs.~(\ref{eqn:defABalpha}), (\ref{eqn:defABalpha2}), (\ref{eqn:Exterf2}),
and (\ref{eqn:Emif2}), 
we arrive at the following definition of the scattering strength.

\begin{definition} The \emph{scattering strength} is a coefficient associated 
with each scatterer that characterizes the mapping from 
\(u_{\alpha}(\bm{r}_{\alpha})\) to \(u_{\alpha}^{\uparrow}(\bm{r})\) 
through the following Green's function: 
\begin{equation}
 u_{\alpha}^{\uparrow} (\bm{r}) = 
 s_{\alpha}(a)\, G_k(\bm{r},\bm{r}_{\alpha})\, u_{\alpha}(\bm{r}_{\alpha}).
 \end{equation}
 Here, \(s_{\alpha}\) is referred to as the \emph{scattering strength}.
  \end{definition}

The explicit forms of the scattering strength can be derived from 
Eqs.~(\ref{eqn:defABalpha}) and (\ref{eqn:defABalpha2}) as
\begin{equation}
s_{\alpha} (a)
= \left \{
\begin{array}{ll}
4\pi a^2  \mu_{\alpha}/3 & \text{in case 1},\\[0.3em]
\displaystyle 4\pi a 
\frac{\mathcal{M}_{\alpha}}{\mathcal{N}_{\alpha} }
& \text{in cases 2a, and 3} \\
\displaystyle -\frac{4\pi i}{k} \Bigl  (1+\frac{a \mathcal{C}_0}{j_0 (\sqrt{\mu_{\alpha} })}\Bigr )
& \text{in  case 2b }, \\ [0.6em]
\end{array} \right .
\label{eqn:strngth1-3}
\end{equation}

One of the important properties of the coefficient 
\(B_{\alpha}\) is that
\begin{equation}
B_{\alpha} \to 0, \quad (a\to 0), \quad \text{(in cases 1, 2a, and 3)}. 
\end{equation}
Specifically, we need to be aware that in case 3, 
\begin{equation}
\mathcal{N}_{\alpha} \ne 0,\ (\forall a >0 ).
\end{equation}
Therefore, we have the
following lemma.

\begin{lemma}
The scattering strength \(s_{\alpha} (a)\) satisfies
\begin{equation}
s_{\alpha} (a) \to \left \{
\begin{array}{ll}
0 & \quad \text{(in cases 1, 2a, and 3)},\\
\displaystyle -\frac{4\pi i}{k} & \quad \text{(in case 2b)},\\
\end{array} \right .
\label{eqn:sascale}
\end{equation}
when \(a \to 0 \).
Therefore, the scattered field vanishes in the zero-size limit for 
cases 1, 2a, and 3, but does not vanish for case 2b.
\end{lemma}

\noindent
\textbf{Remark 2.}
The solution of this single-scattering problem in terms of the Lippmann-Schwinger equation is
\begin{equation}
u_{\alpha}^{\uparrow} (\bm{r} )=
\frac{\mu_{\alpha}}{a^d} 
\int_{V_a (\bm{r}_{\alpha})}
G_k (\bm{r}, \bm{r}')\, u_{\alpha}^{\downarrow}
( \bm{r}' )\, d\bm{r}',
\label{eqn:LS}
\end{equation}
which suggests that the scattering strength is proportional to
\(\mu_{\alpha}/a^d\). 
However, the scattering strength obtained in Eq.~(\ref{eqn:strngth1-3}) is significantly different
and we never arrive at Lemma 2 from Eq. (\ref{eqn:LS}).
The scattering strength is uniquely determined by the transmission and emission 
coefficients of the scatterer and therefore is not, in general, identical to the bare coupling constant.

\subsection{Asymptotic Interaction Equation for \(N\) Scatterers}  
Let us return to Eq.~(\ref{eqn:ub_aE1}), 
which expresses the continuity of the wavefield across the boundary
of the \(\alpha\)-th scatterer. 
Equation~(\ref{eqn:ub_aE1})  can be modified as follows:
\begin{equation}
u_{\alpha} ( \bm{r} )=u_{\alpha}^{\downarrow} (\bm{r} ) -u_{\alpha}^{\uparrow} 
(\bm{r} ),\ \  (\bm{r} \in C_{\alpha} )
\label{eqn:ub_aE2}
\end{equation}
where  \(C_{\alpha} \) is the boundary of the \(\alpha\)-th scatterer 
defined by
\begin{equation}
C_{\alpha} = \{ \bm{r}\in \mathbb{R}^3\ \bigl  | \ 
| \bm{r} - \bm{r}_{\alpha} | =a \bigr . \}.
\end{equation}

We construct the asymptotic interaction equation for the multiple-scattering problem
from Eq. (\ref{eqn:ub_aE2}) 
as \(\bm{r} \to \bm{r}_{\alpha}, \ \bm{r}\in C_{\alpha} \).  
Here,  the left-hand side of Eq.~(\ref{eqn:ub_aE2}) describes
the external field propagating to the boundary of the \(\alpha\)-th scatterer
(i.e., the off-diagonal term), 
and the right-hand side represents the diagonal term
(i.e., self-interaction term) in the interaction equation. 
By using the emission and transmission coefficients, we have the 
following theorem.

\begin{theorem}
The self-interaction term of the interaction equation (i.e., the right-hand side
of Eq. (\ref{eqn:ub_aE2})) 
does not diverge; rather, it is asymptotically equivalent  to the external field
at the center of the \(\alpha\)-th scatterer as
 \(\bm{r} \to \bm{r}_{\alpha}, \ \bm{r}\in C_{\alpha} \), as long as 
 the transmission and emission coefficients are employed 
 for its evaluation, which verifies the consistency of Eq.~(\ref{eqn:ub_aE2}). 
\end{theorem}

\begin{proof}
According to Eqs.~(\ref{eqn:Emif1}) and (\ref{eqn:Intf1}), 
the right-hand side of Eq. (\ref{eqn:ub_aE2}) can be  modified as follows: 
\begin{equation}
u_{\alpha}^{\downarrow}( \bm{r} )- u_{\alpha}^{\uparrow} (\bm{r} )
=\frac{1}{\sqrt{4\pi}}j_{0} (k_{\alpha} a ) A_{\alpha} 
- \frac{1}{\sqrt{4\pi}}h_{0}^{(2)} (k a ) B_{\alpha}. 
\label{eqn:diag2}
\end{equation}
Substituting  the explicit forms of  the transmission and emission
coefficients, 
\(A_{\alpha}\) and \(B_{\alpha} \) given
in Eqs. (\ref{eqn:defABalpha}) and (\ref{eqn:defABalpha2}),  into 
Eq. (\ref{eqn:diag2}) yields
the following:
\begin{eqnarray}
\lefteqn{ 
u_{\alpha}^{\downarrow} (\bm{r} ) -u_{\alpha}^{\uparrow} 
(\bm{r} )}
\nonumber \\
&=& \left \{
\begin{array}{ll}
u_{\alpha} ( \bm{r}_{\alpha} ) + O \bigl ((ka) \bigr )
&  \hspace * {3em} {\text{in cases\ 1,\ 2a,\  and\ 2b} } \\
u_{\alpha} (  \bm{r}_{\alpha} ) + O \bigl ( (ka)^2 \bigr )
&  \hspace * {3em} {\text{in case\ 3} }
\end{array} \right .
\nonumber \\
& & 
\label{eqn:daig}
\end{eqnarray}
To derive Eq.~(\ref{eqn:daig}), we use 
Eqs.~(\ref{eqn:defABalpha}) and (\ref{eqn:defABalpha2}) and
the following asymptotic approximation:
\begin{equation}
h_{0}^{(2)} (ka)=
\left \{
 \begin{array}{ll}
 \displaystyle \frac{i}{ka} + O(1) & \hspace * {3em} {\text{ Cases\ 1 and \ 2a, } }\\
 & \\
 \displaystyle  \frac{i}{ka} + 1 + O( (ka)) & \hspace *
{3em} {\text{ Cases\ 2b\  and\ 3 } }\\
\end{array} \right .
\end{equation}
\end{proof}

The left-hand  side of Eq. (\ref{eqn:ub_aE2}), which is the external field 
propagating to the boundary of the \(\alpha\)-th scatterer,    
can be represented by
\begin{equation}
  u_{\alpha} (\bm{r} )=
  \sum_{\beta=1, \beta \ne \alpha}^{N}
  s_{\beta} (a ) G_k (\bm{r}, \bm{r}_{\beta} )
  u_{\beta} (\bm{r}_{\beta} ) + G_k(\bm{r}, \bm{r}_{F} )
,\ (\bm{r} \in C_{\alpha} )
\label{eqn:Extf3G}
\end{equation}
where \(\bm{r}_{F}\) is the position of a point force that generates
the incident wavefield.

The asymptotic behavior of the Green's function 
\(\bm{r} \sim \bm{r}_{\alpha} \) in Eq. (\ref{eqn:Extf3G}) is
as follows:
\begin{eqnarray}
G_k( \bm{r}, \bm{r}_{\beta} )
&=&
G_k( \bm{r} -\bm{r}_{\alpha}, \bm{r}_{\beta} - \bm{r}_{\alpha} )
\nonumber \\
&=& -ik  \sum_{n=0}^{\infty}
\sum_{|m| \le n} j_n ( k a ) Y_{n}^{m} ( \theta, \varphi )
h_{n}^{(2)} (k | \bm{r}_{\beta} - \bm{r}_{\alpha} | )
Y_{n}^{m} ( \theta ' , -\varphi')
\nonumber \\
&\simeq &  -\frac{ik}{4 \pi} \: j_0 (ka ) h_{0}^{(2)} (k | \bm{r}_{\beta} - \bm{r}_{\alpha} | ), 
\quad (n=0 \: \mbox{approximation})
\nonumber \\
&=& j_0 (ka) G_k(\bm{r}_{\alpha} , \bm{r}_{\beta} ) 
\nonumber \\
&=&
G_k(\bm{r}_{\alpha} , \bm{r}_{\beta} ) + O\bigl ((ka)^2 \bigr ), 
\nonumber \\
& &  \Bigl (\because \ j_0(ka)=1+O\bigl((ka)^2\bigr ) \Bigr )
\label{eqn:approgr}
\end{eqnarray}
Likewise, 
\begin{eqnarray}
G_k( \bm{r}, \bm{r}_{F} )
&=&
G_k( \bm{r}_{\alpha}, \bm{r}_{F} )+ O\bigl ((ka)^2 \bigr ), \ \  (\bm{r}\in C_{\alpha} )
 \end{eqnarray}
 
Therefore, by omitting the \( (ka)\) term for cases 1, 2a, and 2b, or the \( (ka)^2 \) term for
case 3, we have the following corollary.
\begin{corollary}
For cases 1, 2a, 2b, and 3, the Foldy algebraic equation emerges
as the leading-order interaction equation
in the following form:
\begin{equation}
\sum_{\beta=1, \beta \ne \alpha}^{N}
  s_{\beta} (a ) G_k (\bm{r}_{\alpha}, \bm{r}_{\beta} )
  u_{\beta} (\bm{r}_{\beta} ) + G_k(\bm{r}_{\alpha}, \bm{r}_{F} )
=u_{\alpha} (\bm{r}_{\alpha} ).
\label{eqn:FoldyPr}
\end{equation}
 \end{corollary} 
 
The total wavefield outside the \(N\) point-like scatterers is given by the sum of
the scattered waves emitted from all scatterers and the incident wavefield.
Therefore, the total wavefield can be represented  by
\begin{equation}
u ( \bm{r} )= \sum_{\beta=1}^{N}
  s_{\beta} (a ) G_k (\bm{r}, \bm{r}_{\beta} )
  u_{\beta} (\bm{r}_{\beta} ) + G_k(\bm{r}, \bm{r}_{F} ), \ \ 
  ( \bm{r}\in \mathbb{R}^3\setminus Y )
\label{eqn:totf}
\end{equation}
where
\begin{equation}
Y =\{ \bm{r}_1, \bm{r}_2, \ldots , \bm{r}_N \}
\label{eqn:setY}.
\end{equation}
Substituting Eq.~(\ref{eqn:FoldyPr}) into Eq.~(\ref{eqn:totf}) yields
\begin{eqnarray}
u ( \bm{r} ) & = & \sum_{\alpha=1}^{N}\sum_{\beta=1}^{N}
\Bigl [ \Gamma_{a,d,Y}^{(F)} (k) \Bigr ]_{\alpha \beta}^{-1}   
   s_{\alpha} (a ) G_k (\bm{r}, \bm{r}_{\alpha})
  G_k (\bm{r}_{\beta}, \bm{r}_{F}  ) 
\nonumber \\
& & 
+ G_k(\bm{r}, \bm{r}_{F} ), \ \ 
( \bm{r}\in \mathbb{R}^3\setminus Y )
\label{eqn:totfx}
\end{eqnarray} 
where
\( \Gamma_{a,d,Y}^{(F)} (k) \) is defined by
\begin{equation}
\Bigl [\Gamma_{a,d,Y}^{(F)} (k)\Bigl ]_{\alpha\beta} = \left \{
\begin{array}{ll}
1,  & (\alpha = \beta), \\
-s_{\beta} (a) G_k(\bm{r}_\alpha, \bm{r}_\beta), & (\alpha \ne \beta).
\end{array} \right .
\label{eqn:defGammaX}
\end{equation}
Equation~(\ref{eqn:totfx}) provides a Foldy-type representation of the
solution of Eq.~(\ref{eqn:e1}) obtained within the framework of classical
wave theory.
This representation is valid provided that the matrix
\( \Gamma_{a,d,Y}^{(F)}(k) \)
is invertible.
The invertibility of this matrix is established in the following section. 

 \section{Comparison of Classical 
Foldy Equation with Quantum Mechanics in View of Renormalization Theory}

\subsection{Strong Resolvent Convergence of Foldy Operator
for Cases 1, 2a, and 3}

Recall that our analysis began with Eq.~(\ref{eqn:e1}), 
to which classical wave theory was applied.
Although this equation describes acoustic wave propagation, 
its mathematical treatment is closely related to that of the 
Schr{\"{o}}dinger equation with \( N\)  potential wells,
whose operator is given by
\begin{equation}
H_{a,d,Y} = -\Delta - a^{-d} \sum_{\alpha=1}^{N}
V_{\alpha}\! \bigl (  ( \cdot -\bm{r}_{\alpha} )/a \bigr ),\ \ 
(d=1,2,3)
\label{eqn:Ha}
\end{equation}
where
\begin{equation}
V_{\alpha}\!\bigl (  (\cdot -\bm{r}_{\alpha} )/a \bigr  )
= \mu_{\alpha}
X_{a}(\bm{r} - \bm{r}_{\alpha}).
\label{eqn:am3pot}
\end{equation}

Here, \( \Delta  \) 
is a self-adjoint operator
on the domain \(H^{2,2} ( \mathbb{R}^3)\)
and each \(V_{\alpha} \)
is a real-valued multiplication operator with the property
\begin{equation}
\sup_{\bm{r} \in \mathbb{R}^3}
\bigl| V_{\alpha}(\bm{r} - \bm{r}_{\alpha}) \bigr| < \infty.
\end{equation}
Consequently, the operator defined by Eq. (\ref{eqn:Ha})  
is already self-adjoint on the domain
\(H^{2,2} ( \mathbb{R}^3)\)
(see, for example, Amrein, 1981, p.~46).

This observation shows that, for scatterers of finite spatial extent,
the redefinition  of the coupling constant in the Foldy  equation
is not required to resolve issues of self-adjointness.
Rather,  this redefinition arises from the properly defined transmission and emission
coefficients. 
Importantly, in this setting, there is no need to introduce a family of self-adjoint extensions
of \( H_{a,d,Y}\) parametrized by a real constant.

For the scaling \(a^{-d}\) of the potential, the case \(d=2\) is special. 
The arguments employed by Albeverio et al.~(1988), summarized in the Appendix , 
can be directly applied to the limit process \(a\to 0\) for Eq.~(\ref{eqn:Ha}).
Nevertheless, we can construct the resolvent of Eq.~(\ref{eqn:Ha}) 
for the limit process of
\( a\to 0\) using the Foldy  equation approach for every case. 
 
In this section, we discuss the resolvent kernel for cases 1, 2a, and 3,
which converges to the resolvent kernel of the  free Hamiltonian as \(a\to 0\). 
According  to Eq.~(\ref{eqn:totfx}), 
the resolvent kernel of Eq.~(\ref{eqn:Ha})  for  cases 1,2a, and 3,  can be expressed as
\begin{eqnarray}
\lefteqn{
\Bigl(  H_{a,d,Y}   - k^2 \Bigr)^{-1} (\bm{r}, \bm{r}')
}
\nonumber \\
&=&
G_k(\bm{r}, \bm{r}') +
\sum_{\alpha=1}^{N}\sum_{\beta=1}^{N}
\Bigl [ \Gamma_{a,d,Y}^{(F)} (k) \Bigr ]_{\alpha \beta}^{-1} 
s_{\alpha}  (a) G_k ( \bm{r}, \bm{r}_{\alpha} ) 
G_k (\bm{r}_{\beta}, \bm{r}' ), \nonumber  \\
& &
\bigl ( k^2 \in \rho ( H_{a,d,Y}), \ {\rm{Im}}\: (k)<0  \bigr ).
\label{eqn:resFol}
\end{eqnarray}
Since Eq.~(\ref{eqn:totfx}) represents the wavefield generated by a
point source located at \(\bm{r}_F\), the corresponding resolvent
kernel is obtained by replacing the source point \(\bm{r}_F\) with a
general variable \(\bm{r}'\), which yields Eq.~(\ref{eqn:resFol}).

An important property of the matrix \( \Gamma_{a,d,Y}^{(F)} (k) \) 
in Eq.~ (\ref{eqn:resFol})
is summarized in the following lemma.\\

\begin{lemma}
The matrix \( \Gamma_{a,d,Y}^{(F)}(k) \) 
for cases 1, 2a, and 3
is invertible and its inverse is 
uniformly bounded with respect to the scale parameter \( a \).
\end{lemma}

\begin{proof}
Define the matrix
\begin{equation}
D_a := I - \Gamma_{a,d,Y}^{(F)}(k),
\end{equation}
where \( I \) denotes the identity matrix. By Corollary~1, 
the scattering strengths satisfy \( s_\alpha(a) \to 0 \) as \( a \to 0 \), and hence
\begin{equation}
\| D_a \| = \| I - \Gamma_{a,Y}^{(F)}(k) \| \to 0 \quad  (a \to 0).
\end{equation}
Therefore, for sufficiently small \( a \), we have \( \| D_a \| < 1 \), 
which implies that \( \Gamma_{a,Y}^{(F)}(k) \) is invertible via Neumann series expansion.
Moreover,
there exists a constant \(C>0\), independent of \(a\), such that 
\[ \| \Gamma_{a,d,Y}^{(F)}(k)^{-1} \| \le C . \]
This proves that \( \Gamma_{a,d,Y}^{(F)}(k)^{-1} \) is uniformly bounded 
with respect to the scale parameter \( a \).
\end{proof}

As a result, we have the following theorem that characterizes
the convergence of the Foldy resolvent.

\begin{theorem}
The Foldy resolvent represented by Eq.~(\ref{eqn:resFol}) converges in the strong resolvent sense:
\begin{eqnarray} 
\lefteqn{
\mathop{\mathrm{s\!-\!lim}}_{a\to 0}
\left( H_{a,d,Y} -k^2 \right )^{-1} 
}\nonumber \\
& = & \left( -\Delta  - k^2 \right )^{-1},\ \ 
 \left( k^2 \in \rho ( - \Delta ),\ \mathrm{Im}\,(k) < 0 \right),
\end{eqnarray} 
for cases 1, 2a, and 3.
In other words, in the sense of strong resolvent convergence, 
the interaction part vanishes as 
\( a \to 0 \).
\end{theorem}

\begin{proof}
Let \( f \in L_2(\mathbb{R}^3) \) and define \( g \) as
\begin{eqnarray}
g(\bm{r}) &=& \left( \left( H_{a,d,Y} - k^2 \right)^{-1} f - 
\left( -\Delta - k^2 \right)^{-1} f \right)(\bm{r}) \nonumber \\
&=& \sum_{\alpha=1}^{N} \sum_{\beta=1}^{N}
\left[ \Gamma_{a,d,Y}^{(F)}(k) \right]_{\alpha \beta}^{-1} s_{\alpha}(a)
G_k(\bm{r}, \bm{r}_{\alpha}) G_{f,k}(\bm{r}_{\beta}),
\nonumber \\
& & 
\end{eqnarray}
where
\begin{equation}
G_{f,k}(\bm{r}_{\beta}) = \int_{\mathbb{R}^3} G_k(\bm{r}_{\beta}, \bm{r}') f(\bm{r}') 
\, d\bm{r}'.
\end{equation}

By the Cauchy-Schwarz inequality, we have
\begin{equation}
|G_{f,k}(\bm{r}_{\beta})| \le W_{\beta}^{1/2} \|f\|_{L_2(\mathbb{R}^3)},
\end{equation}
where
\begin{equation}
W_{\beta} = \int_{\mathbb{R}^3} |G_k(\bm{r}_{\beta}, \bm{r}')|^2 \, d\bm{r}'.
\label{eqn:defGfk}
\end{equation}

According to the explicit form of Green's function,
\begin{equation}
|G_k(\bm{r}_{\beta}, \bm{r}')|^2 = \frac{\exp(-2k_I R)}{16\pi^2 R^2}, \quad (k_I > 0),
\label{eqn:Gk2}
\end{equation}
where \( R = |\bm{r}_{\beta} - \bm{r}'| \) 
and \( k = k_R - i k_I \). This shows that \( W_{\beta} \) is finite and independent of \( a \).

Therefore, the \( L_2 \)-norm of \( g \) is estimated as
\begin{eqnarray}
\|g\|_{L_2(\mathbb{R}^3)}^2
&=& \left\| \left( H_{a,d,Y} - k^2 \right)^{-1} f -
 \left( -\Delta - k^2 \right)^{-1} f \right\|_{L_2(\mathbb{R}^3)}^2 \nonumber \\
&\le& \left( \sum_{\alpha=1}^{N} |s_{\alpha}(a)|^2 W_{\alpha} \right)
\left\| \left[ \Gamma_{a,Y}^{(F)}(k) \right]^{-1} \right\|^2
\nonumber \\
& & \times 
\left( \sum_{\beta=1}^{N} W_{\beta} \right) \|f\|_{L_2(\mathbb{R}^3)}^2 \nonumber \\
&\le& C a^{m(d)} \: (1 + \epsilon)^2 \|f\|_{L_2(\mathbb{R}^3)}^2 \to 0 \quad (a \to 0),
\end{eqnarray}
where
\begin{equation}
m(d)=
\begin{cases}
4, & (d=1),\\
2, & (d=2, 3).
\end{cases}
\end{equation}

This completes the proof.
\end{proof}

\noindent
\textbf{Remark 3.} \\
Within the framework of renormalization theory, which is presented in 
Eqs.~(\ref{eqn:resolfx}) and (\ref{eqn:renorm2}) in the Appendix, 
the result of Theorem~2 may be interpreted as the limiting case 
\(\theta_\alpha=+\infty\). 
 \\

\subsection{Strong Resolvent Convergence of Foldy Operator for Case 2b} 
For case 2b, 
we can use Eq.~(\ref{eqn:resFol}) to represent the resolvent kernel. 
For this case, the scattering strength given by Eq.~(\ref{eqn:strngth1-3}) yields the 
following representation of the resolvent kernel:
\begin{eqnarray} 
\lefteqn{ \Bigl( H_{a,2',Y} - k^2 \Bigr)^{-1} (\bm{r}, \bm{r}') } \nonumber \\ 
&=& G_k(\bm{r}, \bm{r}') + \sum_{\alpha=1}^{N}\sum_{\beta=1}^{N}
 \Bigl [ \Gamma_{a,2',Y}(k) \Bigr ]_{\alpha \beta}^{-1}
 G_k(\bm{r}, \bm{r}_{\alpha}) G_k(\bm{r}_{\beta}, \bm{r}'), \nonumber \\ 
 && \bigl( k^2 \in \rho(H_{a,2',Y}), \ {\rm Im}(k)<0 \bigr), 
 \label{eqn:Ha2y}
 \end{eqnarray}
 where the subscript \(2'\) indicates that the corresponding operator or matrix is
  associated with case~2b.
  Additionally, the matrix \(\Gamma_{a,2',Y}\) is given by 
   \begin{equation} \Bigl [ \Gamma_{a,2',Y}(k) \Bigr ]_{\alpha\beta} =
    \left\{ \begin{array}{ll} \displaystyle \frac{ik}{4\pi} \Bigl( 1+\frac{a\mathcal{C}_0} 
    {j_0(\sqrt{\mu_{\alpha}})} \Bigr), & (\alpha=\beta), \\
    [0.3cm] -G_k(\bm{r}_{\alpha},\bm{r}_{\beta}), & (\alpha\neq\beta). \end{array} \right. 
       \label{eqn:GHay2}  
       \end{equation}
 Notably, Eqs.~(\ref{eqn:Ha2y}) and (\ref{eqn:GHay2}) 
 match the form of Eqs.~(\ref{eqn:resolfx})
 and (\ref{eqn:renorm2}) in the Appendix, respectively. 
 Therefore, the classical wave analysis leads naturally to the renormalization framework.
 
Using Eq.~(\ref{eqn:a9}) in the Appendix, 
we can now establish the strong resolvent convergence of Eq.~(\ref{eqn:Ha2y}) as \(a\to 0\).
\begin{theorem}
The Foldy resolvent defined by Eq.~(\ref{eqn:Ha2y}) converges in the strong resolvent sense:
\begin{eqnarray} 
\lefteqn{
\mathop{\mathrm{s\!-\!lim}}_{a \to 0}
\bigl( H_{a,2',Y} - k^2 \bigr)^{-1}
} \nonumber \\
& =&
\bigl( H_{2',Y} - k^2 \bigr)^{-1},\ 
\bigl ( k^2 \in \rho ( H_{2',Y}), \ {\rm{Im}}\: (k)<0  \bigr )
\end{eqnarray}
whose kernel is given by
\begin{eqnarray} 
\lefteqn{
\bigl( H_{2',Y} - k^2 \bigr)^{-1} (\bm{r}, \bm{r}')
}\nonumber \\
&=& 
G_k(\bm{r}, \bm{r}') +
\sum_{\alpha=1}^{N}\sum_{\beta=1}^{N}
\Bigl [ \Gamma_{Y} (k) \Bigr ]_{\alpha \beta}^{-1} 
G_k ( \bm{r}, \bm{r}_{\alpha} ) 
G_k (\bm{r}_{\beta}, \bm{r}' ), 
 \nonumber  \\
 & &
\label{eqn:resFol3}
\end{eqnarray} 
where
   \begin{equation} 
   \Bigl [ \Gamma_{Y}(k) \Bigr ]_{\alpha\beta} =
    \left\{ \begin{array}{ll} \displaystyle \frac{ik}{4\pi} 
    , & (\alpha=\beta), \\
    [0.3cm] -G_k(\bm{r}_{\alpha},\bm{r}_{\beta}), & (\alpha\neq\beta). \end{array} \right. 
       \label{eqn:GHay3}  
       \end{equation}
 In contrast to cases 1, 2a, and 3, the interaction part does not vanish in the limit \(a\to0\).    
Thus, the limiting Foldy operator corresponds precisely to the renormalized point-interaction 
Hamiltonian with \(\theta_{\alpha}=0\).
\end{theorem}

\noindent
\textbf {Remark 4.}\\
Albeverio et al. (1988) showed that there exist specific cases in which
the interaction part does not vanish in the limit \(a\to0\) when the
potential \(V_\alpha\) belongs to the Rollnik class and \(d=2\).
As mentioned above, their analysis was restricted to the case \(d=2\),
which is associated with the unitary scaling group.
These cases are determined by the properties of the potentials. 
A detailed discussion is given in the next section. 
  
In contrast, within the classical wave-theoretic framework considered in
the present study, we found that a corresponding phenomenon arises when
\(d=2\). In this case, the leading-order asymptotic evaluation of the
transmission and emission coefficients breaks down, and the next-order
term in the asymptotic expansion of Eq.~(\ref{eqn:sol_AB}) must be taken
into account. The inclusion of this second-order contribution leads to a
case in which the interaction part remains nonvanishing in the limit
\(a\to0\). This degeneracy occurs when
\(\cos(\sqrt{\mu_\alpha})=0\), indicating that the corresponding
nontrivial limit is determined by the bare coupling constants.

\section{Discussion}

In the monograph of Albeverio et al. (1988), the authors considered
Hamiltonians of the form
\begin{equation}
H_{a,Y}
=
-\Delta
+
a^{-2}
\sum_{\alpha=1}^{N}
\lambda_{\alpha}(a)
V_{\alpha}
\bigl((\cdot-\bm r_{\alpha})/a \bigr),
\end{equation}
where \(V_{\alpha}\) belongs to the Rollnik class and
\(\lambda_{\alpha}(a)\) is assumed to be real analytic in a
neighborhood of the origin with
\begin{equation}
\lambda_{\alpha}(0)=1,
\qquad
\alpha=1,\ldots,N.
\end{equation}

They showed that the resolvent of \(H_{a,Y}\)
can be represented in the form
\begin{eqnarray}
(H_{a,Y}-k^2)^{-1}
&=&
G_k
-
a
\sum_{\alpha,\beta=1}^{N}
A_{a,\alpha}(k)
\bigl[1+B_a(k)\bigr]^{-1}_{\alpha\beta}
C_{a,\beta}(k),
\nonumber\\
&&
k^2\in \rho(H_{a,Y}),
\quad
{\rm Im}(k)<0,
\end{eqnarray}
where \(A_{a,\alpha}(k)\),
\(B_a(k)\), and
\(C_{a,\beta}(k)\) are Hilbert--Schmidt operators.
Their explicit kernel representations are given in
Albeverio et al. (1988).

The limiting behavior of
\(a[1+B_a(k)]^{-1}\)
is highly nontrivial.
As a consequence, depending on the properties of the potentials,
the limiting resolvent may converge either to that of the free
Hamiltonian or to that of a Hamiltonian with point interactions.
Because of the role of the unitary scaling group,
Albeverio et al. (1988) restricted their analysis to the case \(d=2\).

Although the limiting Hamiltonian in renormalization theory is
described by a family of self-adjoint extensions parameterized by
\( \theta_{\alpha}\), the parameter is not arbitrary when one starts from a given
family of less singular Hamiltonians. 
The limiting procedure selects
a specific value of \(\theta_{\alpha} \)  determined by the properties of the
approximating potentials.

In the present study, we examined the cases
\(d=1,2,3\) from the viewpoint of classical wave theory.
For cases 1, 2a, and 3, the leading-order asymptotic evaluation
of the transmission and emission coefficients remains valid.
Consequently, the scattering strength satisfies
\[
s_{\alpha}(a)\to 0,
\qquad
a\to0,
\]
which leads to strong resolvent convergence to the free Hamiltonian.

As discussed in Remark 4,
there exists a special case in which the interaction part remains
nonvanishing in the zero-size limit.
This occurs in case 2b.
From this viewpoint, the present analysis also identifies
\(d=2\) as the distinguished  scaling, in agreement with the
renormalization framework developed by Albeverio et al. (1988).
A similar phenomenon to renormalization theory
is observed in the present classical wave theory
for selecting the specific 
value of \(\theta_{\alpha}\).

The finite-size scatterer model considered in this work is already
self-adjoint for every finite value of \(a\),
and the limiting procedure determines the parameter of
the limiting Foldy operator through the asymptotic behavior of the
transmission and emission coefficients.
The limiting Foldy operator corresponds precisely to the renormalized point-interaction 
Hamiltonian with \(\theta_{\alpha}=\infty\) for cases 1, 2a, and 3 and with \( \theta_{\alpha} =0\)
for case 2b. 

In the present framework, the Foldy operator is not introduced as an
independent approximation model. Rather, it emerges naturally from the
asymptotic equivalence between the diagonal self-interaction term and
the external field evaluated at the scatterer center. 
Consequently, the Foldy equation and the associated Foldy resolvent
should be understood as asymptotic representations of the original
finite-size scatterer problem.

\section{Conclusion} 
In this study, we compared a classical wave theory  for very small spherical 
scatterers with the renormalization technique 
developed in quantum mechanics. The governing equation considered in this work 
was the wave equation for sound propagation under the assumption that the 
gradient of the mass density is negligible. 
This equation can also be interpreted 
as a Schr{\"o}dinger equation with compactly supported potential wells, 
which enabled us to compare these two approaches. 
The  classical wave theory employed in this article  
was based on spherical harmonics expansions 
and asymptotic analysis. We showed that the monopole approximation is valid 
when the spatial scale of the scatterers is sufficiently small. 
One of the significant  findings obtained from the classical wave-theoretic  approach 
was that the diagonal self-interaction contribution
can be represented asymptotically by the external field evaluated at the 
scatterer center as the spatial scale
of the scatterers tends to zero, provided that the transmission and emission 
coefficients are properly defined. 
As a result, the Foldy equation emerges as the leading-order interaction
equation, provided that the bare coupling 
constants are properly redefined via the transmission and emission coefficients. 
Furthermore, the Foldy equation enabled the construction of the resolvent. 
We showed that, in most cases, the resolvent converges 
(in the strong resolvent sense) to that of the free Hamiltonian 
as the spatial scale of the scatterers tends to zero. 
On the other hand, we also identified a specific case in which the interaction 
effect survives in the zero-size limit. At first sight, the classical wave-theoretic approach 
appears to be quite different from the renormalization technique developed 
in quantum mechanics. However, from the viewpoint 
of strong resolvent convergence, the classical wave-theoretic approach can be 
understood within the framework of renormalization theory.

Although the limiting Hamiltonian in renormalization theory
is formally described by a family of self-adjoint extensions,
the limiting procedure selects a specific value of  \(\theta_\alpha\) when one starts
from a given family of less singular Hamiltonians.
A similar selection mechanism occurs in the
classical wave-theoretic framework.
This common limiting structure provides a direct connection
between classical wave theory and renormalization theory.

\section * {Acknowledgments}
This work was supported by the Japan Society for the Promotion of Science KAKENHI under grant number 25K15150.

\renewcommand{\theequation}{A.\arabic{equation}}
\renewcommand{\thesection}{A.\arabic{section}}
\renewcommand{\thefigure}{A.\arabic{figure}}
\setcounter{section}{0}
\setcounter{equation}{0}
\setcounter{figure}{0}

\section*{APPENDIX\\ Overview of Renormalization Theory}

Renormalization theory in quantum mechanics addresses operators of the form:
\begin{equation}
H = -\Delta - \sum_{\alpha=1}^{N} \mu_{\alpha} \delta (\cdot - \bm{r}_{\alpha}),
\label{eqn:renorm1}
\end{equation}
where \( \mu_{\alpha} \) for \( \alpha = 1, \ldots, N \) are referred to as coupling constants. 
For the operator given in Eq.~(\ref{eqn:renorm1}),
the resolvent \( (H-k^2)^{-1} \) is not well-defined, and \( H \) fails to be self-adjoint.
In quantum mechanics, however, self-adjointness of operators is a fundamental requirement, 
as only self-adjoint operators can represent physical observables.

To resolve this difficulty, renormalization theory for Eq.~(\ref{eqn:renorm1})
 (see Albeverio et al., 1988) introduces the 
following closed symmetric operator:
\begin{equation}
\dot{H}_Y = -\Delta \bigl|_{\mathscr{D}_Y},
\end{equation}
where the domain \( \mathscr{D}_Y \) is defined as
\begin{equation}
\mathscr{D}_Y = \left\{
\phi \in H^{2,2}(\mathbb{R}^3) \ \Big| \ \phi(\bm{r}_\alpha) = 0,\ 
\bm{r}_{\alpha} \in Y
\right\}.
\end{equation}
Recall that the set \( Y \) denotes the collection of singular points 
defined in Eq.~(\ref{eqn:setY}).

This closed symmetric operator has deficiency indices \((N,N)\) 
(see Albeverio et al., 1988, Theorem~1.1.2, p.~113), 
which ensures the existence of a family of self-adjoint extensions.
Let \( -\Delta_{\theta,Y} \) denote one such extension of \( \dot{H}_Y \), 
parametrized by \( \theta \).
The resolvent kernel of \( -\Delta_{\theta,Y} \) can then be expressed as (see Albeverio et al., 1988, Eq. (1.1.33), p.~113)
\begin{eqnarray}
\lefteqn{
\Bigl( -\Delta_{\theta,Y} - k^2 \Bigr)^{-1} (\bm{r}, \bm{r}')
}
\nonumber \\
&=&
G_k(\bm{r}, \bm{r}')
+ \sum_{\alpha=1}^{N} \sum_{\beta=1}^{N}
\bigl[ \Gamma_{\theta,Y}(k) \bigr]_{\alpha\beta}^{-1}
\,
G_k(\bm{r} , \bm{r}_\alpha)
G_k( \bm{r}_\beta, \bm{r}'),
\nonumber \\
&&  
\bigl( k^2 \in \rho (-\Delta_{\theta,Y}),\ \mathrm{Im}(k) < 0 \bigr),
\label{eqn:resolfx}
\end{eqnarray}
where \( \rho(\cdot) \) denotes the resolvent set. 
The matrix \( \Gamma_{\theta,Y}(k) \) is given by
\begin{equation}
\bigl[\Gamma_{\theta,Y}(k)\bigr]_{\alpha\beta} =
\left\{
\begin{array}{ll} 
{\displaystyle \theta_\alpha + \frac{ik}{4\pi}}, & \alpha = \beta, \\[1ex]
-G_k(\bm{r}_\alpha, \bm{r}_\beta), & \alpha \ne \beta.
\end{array}  \right.
\label{eqn:renorm2}
\end{equation}

Each parameter \( \theta_{\alpha} \) is real and may take any value on the extended real line:
\begin{equation}
-\infty < \theta_\alpha \le +\infty.
\label{eqn:A6}
\end{equation}
The term \( \theta_\alpha +{ik}/(4\pi) \) in Eq.~(\ref{eqn:renorm2}) plays 
the role of a renormalized coupling constant. 
In particular, the choice \( \theta_{\alpha} = \infty \) corresponds to the 
free (kinetic energy) Hamiltonian, indicating that no interaction is present at the location of 
the \( \alpha \)-th scatterer.

It is important to understand how the Hamiltonian \( -\Delta_{\theta,Y} \) can be 
approximated by Hamiltonians with more realistic, less singular short-range interactions.
For this purpose, Albeverio et al.\ (1988) introduced the following family of Hamiltonians:
\begin{equation}
H_{a, Y} = -\Delta
- a^{-2}
\sum_{\alpha=1}^{N}
\lambda_{\alpha}(a)\, V_\alpha \Bigl(
( \cdot - \bm{r}_{\alpha} )/a \Bigr),
\qquad (a > 0),
\label{eqn:a7}
\end{equation}
where \( a \) represents the scale parameter for the potentials 
and \( V_{\alpha} \in R \), with \( R \) denoting the Rollnik class (Albeverio et al., 1988, p.~121). 
The coefficients \( \lambda_\alpha(a) \) are assumed to be real analytic in a 
neighborhood of zero and satisfy
\begin{equation}
\lambda_{\alpha}(0) = 1, \qquad (\alpha = 1, \ldots, N).
\end{equation}

Albeverio et al.\ (1988) proved that in the limit \( a \to 0 \), 
this family converges to \( -\Delta_{\theta,Y} \) in the sense that
\begin{eqnarray}
\lefteqn{
\mathop{\mathrm{s\!-\!lim}}_{a \to 0}
\bigl( H_{a,Y} - k^2 \bigr)^{-1}
} \nonumber \\
& = &
 \bigl( -\Delta_{\theta,Y} - k^2 \bigr)^{-1},\ 
\bigl( k^2 \in \rho(-\Delta_{\theta,Y}),\ \mathrm{Im}\,k < 0 \bigr).
\label{eqn:a9}
\end{eqnarray}
Here, the explicit scale dependence of the potential through the factor \( a^{-2} \) is 
crucial for defining the unitary scaling group.

In the main text, we return to the Foldy  equation and discuss its 
renormalization from the perspective of spatial scale dependence in classical wave theory.

\section *
{REFERENCES} 

\begin{list}{}{
\setlength{\leftmargin}{10mm}
\setlength{\listparindent}{-10mm}
\setlength{\itemindent}{-10mm}
\setlength{\labelsep}{-10mm}
\setlength{\labelwidth}{0mm}
}
\item 
Albeverio, S., Gesztesy, F., H{\o}egh-Krohn, R. and
Holden, H. (1988). Solvable models in quantum mechanics,  
Springer-Verlag, (1988).
\item 
Amrein, W.O. (1981). Non-Relativistic Quantum Dynamics, 
Kluwer. 
\item 
Angel, Y.C. and Koba,Y.K., (1998). 
Complex-valuedwavenumber, reflection and transmission
in an elastic solid containing a cracked slab region,
{Int.J. Solids Structures}, 35 573-592
\item 
Arist{\'e}gui,C. and Angel, Y.C. (2002). New results for
isotropic point scatterers: Foldy revisited, 
{Wave motion}, 36, 383-399.
\item %
Challa, D.P. and Sini, M. (2014).
On the justification of the Foldy-Lax approximation
for acoustic scattering by rigid bodies of arbitrary shapes,
{Multiscale Modeling and Simulation}, 12, 55-108
\item 
Foldy, L.L. (1945). 
The multiple scattering of waves. I. General theory
of isotropic scattering by randomly distributed scatterers., 
{Physical  Review}, 67, 107-119. 
\item 
Hu, G. and Sini, M., (2013). 
Elastic scattering by finitely many point-like obstacles,
{Journal of Mathematical Physics}, 54, 042901
\item 
Hu, G., Mantile, A., Sini, M.  and Yin, T. (2020).
Direct and inverse time-harmonic elastic
scattering from point-like and extended obstacles,
{Inverse Problems and Imaging},  1025-1056. doi: 10.3934/ipi.2020054 
\item 
Lax,M. (1951).Multiple scattering of waves. {Reviews of Modern Physics} 23, 287-310.
\item 
Lax, M. (1952).  Multiple scattering of waves II. The effective field in dense systems. 
{Physical Review}, 85, 621-629.
\item 
Lax, M. (1954). Localized perturbations. {Physical Review}, 94, 1391-1392.
\item 
Martin,P.A., (2006). Multiple scattering, Interaction of time-harmonic waves
with $N$ obstacles, Cambridge University Press.
\item 
Ueda,T. (2022). Design concept verification for a transcranial
acoustic lens with instantaneous adaptability, (Journal of  Applied  Physics), 132, 144504.
\item 
Waterman,P.C., Truell,R., (1961).
Multiple scattering of waves, {Journal of Mathematical Physics}, 2,512-537. 
\end{list}

\end{document}